\documentclass[conference]{IEEEtran}
\IEEEoverridecommandlockouts
\usepackage{cite}
\usepackage{amsmath,amssymb,amsfonts}
\usepackage{algorithmic}
\usepackage{graphicx}
\usepackage{textcomp}
\usepackage{xcolor}
\usepackage{mathtools}
\usepackage{siunitx}
\usepackage[textsize=tiny,disable]{todonotes}

\usepackage{layouts}
\usepackage{url}
\def\BibTeX{{\rm B\kern-.05em{\sc i\kern-.025em b}\kern-.08em
    T\kern-.1667em\lower.7ex\hbox{E}\kern-.125emX}}
\newcommand{\transpose}{\textrm{T}}
\newcommand{\defas}{\triangleq}
\newcommand{\rank}[1]{\,\mathrm{rank}\left\{\,#1\,\right\}}

\newcommand{\V}[1]{\ensuremath{\boldsymbol{#1}}}
\newcommand{\M}[1]{\ensuremath{\mathsf{#1}}}

\newcommand{\Vt}[1]{\ensuremath{\tilde{\V{#1}}}}

\newcommand{\Mu}[2]{\ensuremath{\M{#1}_{\text{#2}}}}
\newcommand{\Vh}[1]{\ensuremath{\hat{\V{#1}}}} 
\newcommand{\Vu}[2]{\ensuremath{\V{#1}_{\text{#2}}}}
\newcommand{\Sh}[1]{\ensuremath{\hat{#1}}} 

\newcommand{\Su}[2]{\ensuremath{#1_{\text{#2}}}}

\newcommand{\md}{\,\ensuremath{\mathrm{d}}}
\newcommand{\dx}[2]{\ensuremath{\frac{\md #1}{\md #2}}}

\DeclarePairedDelimiter\norm{\lVert}{\rVert}
\DeclareMathOperator{\cond}{cond}
\DeclareMathOperator{\std}{std}
\newcommand{\Eq}{Eqn.~}
\newcommand{\Fig}{Fig.~}

\newcommand{\Sec}{Sec.~}

\begin{document}
%
\title{Simultaneous Approximation\\of Measurement Values and Derivative Data\\using Discrete Orthogonal Polynomials}
\author{\IEEEauthorblockN{Roland Ritt}
\IEEEauthorblockA{\textit{Chair of Automation} \\
\textit{University of Leoben}\\
Leoben, Austria\\
roland.ritt@unileoben.ac.at}
\and
\IEEEauthorblockN{Matthew Harker}
\IEEEauthorblockA{\textit{Chair of Automation} \\
	\textit{University of Leoben}\\
	Leoben, Austria\\
	matthew.harker@unileoben.ac.at}
\and
\IEEEauthorblockN{Paul O'Leary}
\IEEEauthorblockA{\textit{Chair of Automation} \\
	\textit{University of Leoben}\\
	Leoben, Austria\\
	paul.oleary@unileoben.ac.at}
\and

}
\maketitle
%
\begin{abstract}
This paper presents a novel method for polynomial approximation (Hermite approximation) using the fusion of value and derivative information.
Therefore, the least-squares error in both domains is simultaneously minimized.
A covariance weighting is used to introduce a metric between the value and derivative domain, to handle different noise behaviour.
Based on a recurrence relation with full re-orthogonalization, a weighted polynomial basis function set is generated.
This basis is numerically more stable compared to other algorithms, making it suitable for the approximation of data with high degree polynomials.
With the new method, the fitting problem can be solved using inner products instead of matrix-inverses, yielding a computational more efficient method, e.g., for real-time approximation.

A Monte~Carlo simulation is performed on synthetic data, demonstrating the validity of the method.
Additionally, various tests on the basis function set are presented, showing the improvement on the numerical stability.
\end{abstract}
\section{Motivation}
\begin{figure}[!ht]
	\centering
	\includegraphics[width=\linewidth]{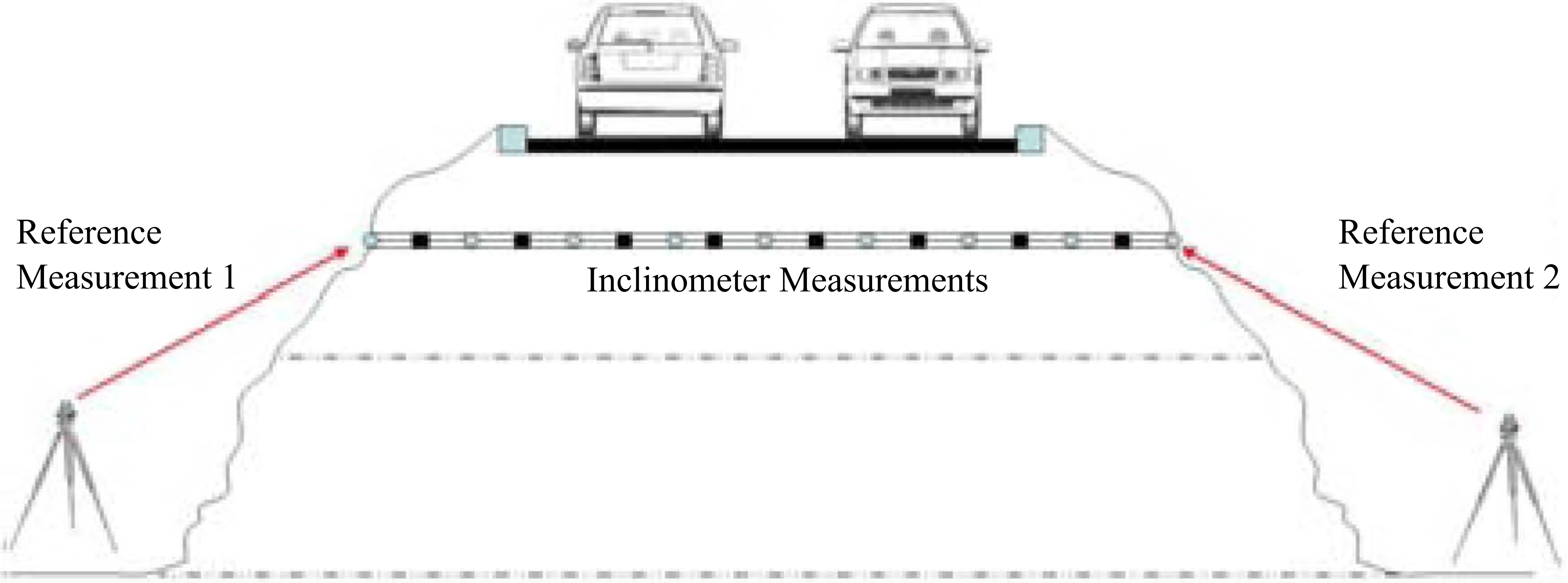}
	\caption{Schematic of geotechnical monitoring using inclination measurements in combination with reference measurements}
	\label{fig:InclinationMeasurement}
\end{figure}

Measurements with a fusion of value and derivative information are common in geotechnical monitoring, e.g. \Fig\ref{fig:InclinationMeasurement}.
In this case, inclinometer measurements and reference measurements from total stations are used to reconstruct the measurement line and detect unwanted changes or deflections. 
These measurements are perturbed by noise  of different characteristics, i.e., the reference measurements act as constraints and are normally more precise than inclinometer measurements.
Especially in engineering problems, polynomials are the model of choice for the approximation of such data, due to their properties and their relation to the underlying physical model.
Within this work, we develop a framework for the approximation with polynomials of given values and derivatives with different noises characteristics using covariance weighted discrete orthogonal polynomials.
The developed method proves to be suitable for high degree polynomial approximation, since it is numerically more stable than other methods.

The main contribution of this paper span:
\begin{enumerate}
	\item The derivation of a novel methodology for the generation of a discrete orthogonal polynomial basis function set which can be used to approximate values and derivatives which are subjected to noise of different levels.
	It uses a three-term recurrence relation with full re-orthogonalization (to increase numerical stability) together with covariance weighting of the residual (for introducing a metric between the value and derivative domain).
	\item The derivation of covariance propagation for the coefficients and the resulting approximation.
	\item The introduction of various measures for the evaluation of the numerical stability.
\end{enumerate}
The paper is structured as follows:
A review of literature is found in \Sec\ref{sec:PreviousWork}.
In \Sec\ref{sec:Framework} all the derivations for the generation of a covariance weighted discrete orthogonal polynomial basis function set are introduced. 
The derivation of the covariance propagation is presented in \Sec\ref{ssec:CovProp}.
To demonstrate the validity of the novel method, a Monte~Carlo experiment for high degree polynomial approximation is performed on a synthetic data set in \Sec\ref{sec:NumExample}.
It is shown that the residual between the data and the approximated values is reduced to noise with the same parameters as used in the generation of the synthetic noise.
The numerical stability of the generated basis is validated in \Sec\ref{ssec:NumQuality} demonstrating the novel method to be advantageous compared to using Vandermonde type basis functions, especially for higher degrees.
Various quality measures are tested and presented for both, complete and incomplete basis function sets.

\section{Review of Literature}
\label{sec:PreviousWork}
The use of value and derivative information for polynomial approximation (Hermite approximation) is not common in literature, whereas the use value and derivative information for interpolation is well-known in literature. 

Based on the idea of a Newton type interpolation, Hermite introduced a similar methodology for the interpolation of values and derivatives.
The original idea can be found in \cite{Hermite1877}.
Based on this, several books, e.g. \cite{Ralston2001}, introduce this idea to a broader audience. 

A generalized version for arbitrary orders of derivatives at given points can be found in \cite{Spitzbart1960}.
This idea is extended to bivariate functions in \cite{Chawla1974}.
An extensive study on multivariate Hermite interpolation can be found in \cite{Sauer1995}.

Since the calculation of those interpolating polynomials is straight forward, Hermite type interpolating polynomials gained popularity for the approximation of functions using a two~point approach rather than Taylor expansion \cite{Shustov2015,Shustov2015a}.
They are well studied in terms of error bounds \cite{Varma1987,Agarwal1991} and compared to other standard methods \cite{Sharifi2011}.
They are used to solving ordinary differential equation \cite{Mennig1983,Grundy2007} and partial differential equation  \cite{Grundy2006,Grundy2005,Grundy2003} by approximating the equations using Hermite type interpolation.
It is pointed out, that using this method lead to higher order of approximation, which improves step-size for iterative methods compared to standard methods.

Another use of the Hermite interpolation is presented in \cite{Komargodski2006}.
In this work it is used to perform a moving-least-squares approximation based on local basis functions which use value and derivative information.

The idea of including derivative information within the approximation of data is mostly found in constrained polynomial approximation, e.g., \cite{Harker2013a,OLeary2014,Gugg2013,OLeary2019}. There it is assumed, that some information (e.g. reference points) is \SI{100}{\percent} certain.
Approximation of data with uncertain value and derivatives is introduced in \cite{Roderick2010}.

The idea of including a covariance weighting on the residuals to get a valid metric between the value and derivative domain is introduced in \cite{OLeary2014} and extended to higher order derivatives in \cite{Ritt2019}.

Within this paper we use this idea of covariance weighting within the generation of a discrete orthogonal polynomial basis function set to improve the numerical stability, making it suitable for high degree polynomial fitting.
\todo{matlab toolbox?}
\section{Theoretical Framework}
\label{sec:Framework}
The herein presented framework addresses the approximation of a polynomial model given noisy values and collocated derivatives using discrete orthogonal polynomials, i.e., systems of polynomials that satisfy a discrete orthogonality constraint.
\subsection{Modelling of Measured Values}
The measured noisy values and derivatives are collected in the vectors  $\Vh{y}=\begin{bmatrix} \Sh{y}_1 & \Sh{y}_2 & \dots & \Sh{y}_n\end{bmatrix}^\transpose$ and $\Vh{y}'=\begin{bmatrix} \Sh{y}_1' & \Sh{y}_2' & \dots & \Sh{y}_n'\end{bmatrix}^\transpose$. They are sampled at the positions  $\V{x}=\begin{bmatrix} x_1 & x_2 & \dots & x_n\end{bmatrix}^\transpose$.
For modelling the measurement, we assume that we measure the true value plus covariant noise. This covariant noise is generated from a vector of gaussian random variables with zero mean and unit variance (i.i.d. noise) together with the according covariance matrix.
Mathematically this is described as 
\begin{IEEEeqnarray}{rCl}
	\label{eq:noisyValueModel}
	\Vh{y} &=& \V{y} + \M{\Lambda}_{\V{y}}^\frac{1}{2} \V{s}\\
	\label{eq:noisyDerivativeModel}
	\Vh{y}' &=& \V{y}' + \M{\Lambda}_{\md\V{y}}^\frac{1}{2} \V{t},
\end{IEEEeqnarray}
where $\V{y}$ and $\V{y}'$ are the true values, $\M{\Lambda}_{\V{y}}$ and $\M{\Lambda}_{\md \V{y}}$ are the associated covariance matrices and $\V{s}$ and  $\V{t}$ are vectors of i.i.d.~noise.
\subsection{Approximation of Values and Derivatives}
The goal is now to approximate the given data in a least-squares sense, using a set of discrete basis functions collected in the columns of the matrix \M{B} and their derivatives \M{B}'.
The true values $\V{y}$ and derivates  $\V{y}'$ are modelled as
\begin{IEEEeqnarray}{rCl}
	\label{eq:valueApprox}
	\V{y} &=& \M{B}\V{\gamma}\\
	\label{eq:derivApprox}
	\V{y}' &=& \M{B}'\V{\gamma},
\end{IEEEeqnarray}
where $\M{B}'$ denotes the first derivative of the basis functions $\M{B}$ with respect to $x$ and $\V{\gamma}$ is the coefficient vector.
To apply Gauss's least-squares theorem, the involved errors must be i.i.d., so we solve for the gaussian random variables  $\V{s}$ and  $\V{t}$ using \Eq\eqref{eq:noisyValueModel}-\eqref{eq:derivApprox}, yielding
\begin{IEEEeqnarray}{rCl}
	\V{s} &=&  \M{\Lambda}_{\V{y}}^ {-\frac{1}{2}} \left(\Vh{y} - \M{B}\V{\gamma}\right)\\
	\V{t} &=&  \M{\Lambda}_{\md \V{y}}^ {-\frac{1}{2}} \left(\Vh{y}' - \M{B}'\V{\gamma}\right).
\end{IEEEeqnarray}
This results in the following functional to be minimized,
\begin{IEEEeqnarray}{rCl}
\label{equ:CostFunction}
E\left(\V{\gamma}\right) &=& \norm{\V{s}}_2^2 + \norm{\V{t}}_2^2 = \\
 &=& \norm{\M{\Lambda}_{\V{y}}^{-\frac{1}{2}}\left(\Vh{y} - \M{B}\V{\gamma}\right)}_2^2 + \norm{\M{\Lambda}_{\md\V{ y}}^{-\frac{1}{2}}\left(\Vh{y}' - \M{B}'\V{\gamma}\right)}_2^2.
\end{IEEEeqnarray}
Clearly, this introduces covariance weighting on the residual and thus a metric between value and derivative domain is established as stated in \cite{OLeary2014} and\cite{Ritt2019}.
Substituting $\M{W}_{\V{y}}^{\frac{1}{2}} \defas \M{\Lambda}_{\V{y}}^{-\frac{1}{2}}$ and $\M{W}_{\md\V{y}}^{\frac{1}{2}} \defas \M{\Lambda}_{\md\V{y}}^{-\frac{1}{2}}$ and minimizing  \Eq\eqref{equ:CostFunction} leads to the normal equations for weighted regression \cite{Golub1996},
\begin{equation}
\label{eq:NormalEqu}	\left(\M{B}^\transpose 	\M{W}_{\V{y}} \M{B}   + \M{B}'^\transpose \M{W}_{\md \V{y}} \M{B}' \right)  \V{\gamma} = 
\M{B}^\transpose \M{W}_{\V{y}}\Vh{y} + \M{B}'^\transpose \M{W}_{\md \V{y}}\Vh{y}'.
\end{equation}
To solve this equation for $\V{\gamma}$ one can use standard methods, e.g, inverting $\left(\M{B}^\transpose \M{W}_{\V{y}} \M{B} + \M{B}'^\transpose \M{W}_{\md \V{y}} \M{B}' \right) $ which is known to be computational costly.

To overcome this problem, we developed a method to find a polynomial basis function set $\M{P}$ and its derivative $\M{P}'$, which fulfil
\begin{equation}
\label{eq:BIdentity}
\M{P}^\transpose 	\M{W}_{\V{y}} \M{P}   + \M{P}'^\transpose \M{W}_{\md \V{y}} \M{P}' = \M{I},
\end{equation}
i.e., a discrete orthogonality condition.
The matrices  $\M{P} = \begin{bmatrix}\Vu{p}{0} & \Vu{p}{1}& \dots &\Vu{p}{i} &\dots &\Vu{p}{d}  \end{bmatrix}$ and $\M{P}' = \begin{bmatrix}\Vu{p}{0}' & \Vu{p}{1}'& \dots &\Vu{p}{i}'  &\dots \Vu{p}{d}'\end{bmatrix}$ are a collection of discrete polynomial basis functions and their derivatives. Each column $\Vu{p}{i}$ represents a polynomial of degree $i$ and   $\Vu{p}{i}'$ is the first derivative of that discrete polynomial.  They are sorted on increasing degrees yielding a basis function set of degree $d$. Respectively a linear combination of those basis functions yield a polynomial of degree $d$  at most.

The coefficients $\V{\gamma}$ for the approximating polynomial can then be easily calculated from \Eq\eqref{eq:NormalEqu} using,
\begin{equation}
\label{eq:alpha}
\V{\gamma} = 
\M{P}^\transpose \M{W}_{\V{y}}\Vh{y} + \M{P}'^\transpose \M{W}_{\md \V{y}}\Vh{y}',
\end{equation}
which are inner products of the basis functions and the covariance weighted measurements. That is, the coefficient $\gamma_i$ for a certain basis function $\Vu{p}{i}$ of degree $i$ can be directly calculated as
\begin{equation}
	\gamma_i = \Vu{p}{i}^\transpose \M{W}_{\V{y}}\Vh{y} + \Vu{p}{i}'^\transpose \M{W}_{\md \V{y}}\Vh{y}',
\end{equation}
which is similar to the calculation of discrete Fourier series.
Thus, the computational efficiency is improved compared to solving the fitting problem using standard algorithms including matrix inverses.
To calculate the approximated values $\Vt{y}$ and  $\Vt{y}'$, we use the estimated parameters $\V{\gamma}$ within the model equation \eqref{eq:valueApprox} and \eqref{eq:derivApprox}, yielding
\begin{IEEEeqnarray}{rCl}
	\label{eq:valueApprox1}
	\Vt{y} &=& \M{P}\V{\gamma}\\
	\label{eq:derivApprox1}
	\Vt{y}' &=& \M{P}'\V{\gamma},
\end{IEEEeqnarray}

\subsection{Synthesis of Weighted Discrete Orthogonal Basis}
\label{ssec:SynthBasis}
In this section, a novel method for the synthesis of a set of weighted discrete orthogonal polynomials fulfilling \Eq\eqref{eq:BIdentity} is presented.

Two important conditions can be directly derived form the identity in \Eq\eqref{eq:BIdentity}.
These are the normal condition
\begin{equation}
\label{eq:BNormalCondition}
\Vu{p}{k+1} ^\transpose  \M{W}_{\V{y}}	\Vu{p}{k+1}   +  \Vu{p}{k+1}'^\transpose \M{W}_{\md \V{y}} \Vu{p}{k+1}'  = 1
\end{equation}
and the orthogonality condition
\begin{equation}
\label{eq:BOrthogonalityCondition}
\Mu{P}{k}^\transpose  \M{W}_{\V{y}}	\Vu{p}{k+1}   + \Mu{P}{k}'^\transpose \M{W}_{\md \V{y}} \Vu{p}{k+1}'  = \V{0},
\end{equation}
which requires, that a basis function of degree $k+1$ is orthogonal to all basis functions of lower degree.
The matrices  $\Mu{P}{k} = \begin{bmatrix}\Vu{p}{0} & \Vu{p}{1}& \dots & \Vu{p}{k} \end{bmatrix}$ and $\Mu{P}{k}' = \begin{bmatrix}\Vu{p}{0}' & \Vu{p}{1}'& \dots &\Vu{p}{k}' \end{bmatrix}$  collect the basis vectors up to degree $k$.

For the generation of a set of suitable polynomial basis functions  $\M{P} $ and their derivatives $\M{P}' $ up to a certain degree $d$, we use a recurrence relation with full re-orthogonalisation (as studied in \cite{OLeary2008Algebraic,Harker2013a}) together with covariance weighting.

Starting from the classical \emph{three-term recurrence relation} from functional analysis for orthogonal polynomials, e.g.~\cite{Gautschi2004a},  also known as Gram-Schmidt process \cite{Golub1996}, a polynomial $\Su{p}{k+1}\left(x\right)$ of degree $k+1$ can be generated from lower degree polynomials using
\begin{equation}
\label{eq:RecurrenceRel}
	\Su{p}{k+1}\left(x\right) = \alpha x \Su{p}{k}\left(x\right) - \beta \Su{p}{k}\left(x\right) - \gamma  \Su{p}{k-1}\left(x\right).
\end{equation}
The derivative of this polynomial with respect to $x$ is calculated as
\begin{equation}
\label{eq:RecurrenceRelDeriv}
	\Su{p}{k+1}'\left(x\right) = \alpha x \Su{p}{k}'\left(x\right) + \alpha x' \Su{p}{k}\left(x\right) - \beta \Su{p}{k}'\left(x\right) - \gamma  \Su{p}{k-1}'\left(x\right).
\end{equation}
This is the definition for the continuous case.

The discrete formulations of \Eq\eqref{eq:RecurrenceRel} and \eqref{eq:RecurrenceRelDeriv} for a given vector $\V{x} = \begin{bmatrix}
x_1 & x_2 & \dots& x_n
\end{bmatrix}^\transpose$ are\footnote{The operation with the symbol $\circ$ denotes the Hadamard product, i.e., the element-wise product.}
\begin{equation}
	\label{eq:ReccurenceVec}
	\Vu{p}{k+1} = \alpha \V{x} \circ \Vu{p}{k} - \beta \Vu{p}{k} - \gamma  \Vu{p}{k-1}
\end{equation}
and
\begin{equation}
	\label{eq:ReccurenceVecDeriv}
	\Vu{p}{k+1}' = \alpha \V{x} \circ \Vu{p}{k}' + \alpha \V{x}' \circ \Vu{p}{k} - \beta \Vu{p}{k}' - \gamma  \Vu{p}{k-1}'.
\end{equation}
The vectors  $\Vu{p}{i} = \begin{bmatrix} \Su{p}{i}\left(x_1\right) &  \Su{p}{i}\left(x_2\right) & \dots &   \Su{p}{i}\left(x_n\right)\end{bmatrix}^\transpose$ and  $\Vu{p}{i}' = \begin{bmatrix} \Su{p}{i}'\left(x_1\right) &  \Su{p}{i}'\left(x_2\right) & \dots &   \Su{p}{i}'\left(x_n\right)\end{bmatrix}^\transpose$ are the polynomial of degree $i$ and its derivative evaluated at the points $\V{x}$.
As shown in \cite{OLeary2009}, the normal three-term recurrence relation is numerically unstable, so we use a complete re-orthogonalization as suggested in \cite{OLeary2008Algebraic}.
Using this improvement, \Eq\eqref{eq:ReccurenceVec} and \eqref{eq:ReccurenceVecDeriv} read as
\begin{equation}
\label{eq:FullReccurenceVec}
\Vu{p}{k+1} = \alpha \V{x} \circ \Vu{p}{k} - \Mu{P}{k}\V{\beta}
\end{equation}
and
\begin{equation}
\label{eq:FullReccurenceVecDeriv}
\Vu{p}{k+1}' = \alpha \V{x} \circ \Vu{p}{k}' + \alpha \V{x}' \circ \Vu{p}{k} -  \Mu{P}{k}'\V{\beta}.
\end{equation}

Using the substitutions
\begin{equation}
\Vu{u}{k} = \V{x} \circ \Vu{p}{k}
\end{equation}
and
\begin{equation}
\label{eq:SubsVecDeriv}
\Vu{v}{k} = \V{x} \circ \Vu{p}{k}' +\V{x}' \circ \Vu{p}{k},
\end{equation}
\Eq\eqref{eq:FullReccurenceVec} and \eqref{eq:FullReccurenceVecDeriv} read as
\begin{equation}
\label{eq:SubsFullReccurenceVec}
\Vu{p}{k+1} = \alpha \Vu{u}{k} - \Mu{P}{k}\V{\beta}
\end{equation}
and
\begin{equation}
\label{eq:SubsFullReccurenceVecDeriv}
\Vu{p}{k+1}' = \alpha \Vu{v}{k} -  \Mu{P}{k}'\V{\beta}.
\end{equation}

Using the orthogonality condition \eqref{eq:BOrthogonalityCondition} together with \eqref{eq:SubsFullReccurenceVec} and \eqref{eq:SubsFullReccurenceVecDeriv} yields
\begin{equation}
\label{eq:PIdentity}
\left(\Mu{P}{k}^\transpose 	\M{W}_{\V{y}} \Mu{P}{k}   + \Mu{P}{k}'^\transpose \M{W}_{\md \V{y}} \Mu{P}{k}'\right) \V{\beta} =  \alpha \left(\Mu{P}{k}^\transpose  \M{W}_{\V{y}}	\Vu{u}{k}   + \Mu{P}{k}'^\transpose \M{W}_{\md \V{y}} \Vu{v}{k}' \right).
\end{equation}
Since, per definition, the previous generated basis functions have to fulfil the orthogonality and normal condition
\begin{equation}
	\Mu{P}{k}^\transpose 	\M{W}_{\V{y}} \Mu{P}{k}   + \Mu{P}{k}'^\transpose \M{W}_{\md \V{y}} \Mu{P}{k}' = \M{I},
\end{equation}
$\V{\beta}$ can be expressed as
\begin{equation}
	\V{\beta} = \alpha \left(\Mu{P}{k}^\transpose  \M{W}_{\V{y}}	\Vu{u}{k}   + \Mu{P}{k}'^\transpose \M{W}_{\md \V{y}} \Vu{v}{k} \right).
\end{equation}
Using this, \Eq\eqref{eq:SubsFullReccurenceVec} and \eqref{eq:SubsFullReccurenceVecDeriv} can be rewritten as
\begin{equation}
\Vu{p}{k+1} = \alpha\left( \Vu{u}{k} - \Mu{P}{k}\Mu{P}{k}^\transpose  \M{W}_{\V{y}}	\Vu{u}{k} - \Mu{P}{k} \Mu{P}{k}'^\transpose \M{W}_{\md \V{y}} \Vu{v}{k}  \right) =  \alpha \Vu{c}{k+1}
\end{equation}
and
\begin{equation}
\Vu{p}{k+1}' =\alpha\left( \Vu{v}{k} - \Mu{P}{k}'\Mu{P}{k}^\transpose  \M{W}_{\V{y}}	\Vu{u}{k} - \Mu{P}{k}' \Mu{P}{k}'^\transpose \M{W}_{\md \V{y}} \Vu{v}{k}  \right) =  \alpha \Vu{c}{k+1}'.
\end{equation}
\todo{$\Vu{c}{k+1}$ use subscript $k+1$ or only $\V{c}$}
$\Vu{c}{k+1}$ and $\Vu{c}{k+1}'$ represent a basis function and (its derivative) which is orthogonal to all previous basis functions but not yet normed.
To fulfil the norm condition in \Eq\eqref{eq:BNormalCondition}, the scaling factor $\alpha$ is calculated as
\begin{equation}
	\alpha = \sqrt{\frac{1}{\Vu{c}{k+1}^\transpose  \M{W}_{\V{y}}\Vu{c}{k+1} +  \Vu{c}{k+1}'^\transpose\M{W}_{\md \V{y}} \Vu{c}{k+1}' }}.
\end{equation}
\todo{change from $\alpha$ to $\alpha_{k+1}$ or make a sidenote??}
This scaling factor and the coefficient vector $\V{\beta}$ have to be calculated for each newly generated basis function.

\begin{figure}[!ht]
	\centering
	\includegraphics{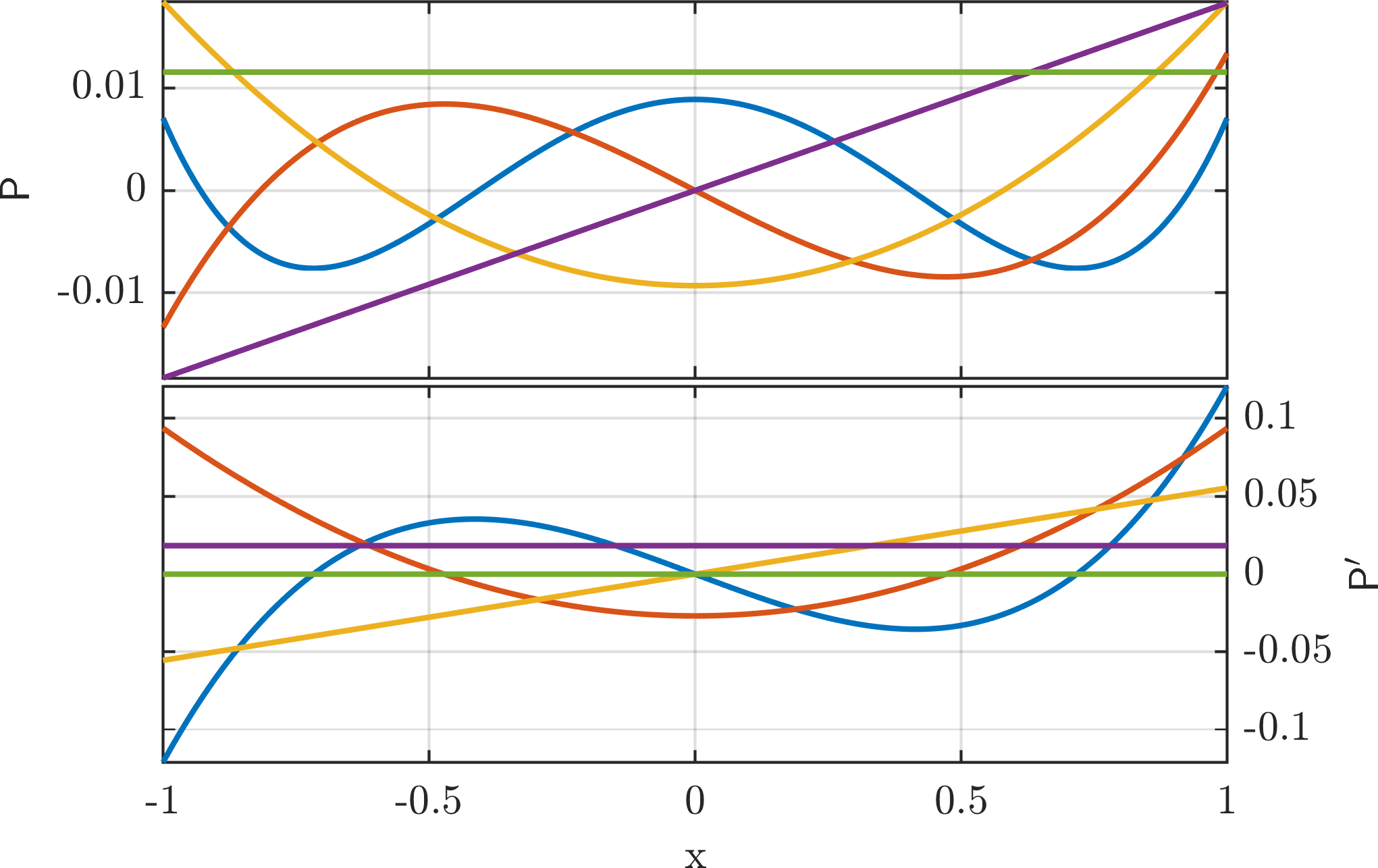}
	\caption{A set of discrete orthogonal basis functions $\M{P}$ of degree $d=4$ and its derivative $\M{P}'$  generated for $n=300$ equally spaced points with $\sigma_{y_i} = \sigma_{y} =  0.2$ and $\sigma_{\md y_i} = \sigma_{\md y} =  0.8$ }
	\label{fig:BasisFns}
\end{figure}

To start the recurrence, the first basis functions $\Vu{p}{0}$, $\Vu{p}{1}$ and their derivatives have to be defined in an initial step, to meet the above mentioned conditions.
The first basis function of degree $d=0$ is defined as
\begin{equation}
\Vu{p}{0} = \frac{\V{e}}{\sqrt{\V{e}^\transpose\M{W}_{\V{y}}\V{e} }},
\end{equation}
which is a normalized constant vector. $\V{e}$ denotes a vector of ones.
The first derivative of this basis function is the zero vector
\begin{equation}
\Vu{p}{0}' = \V{0}.
\end{equation}
For the second basis function $\Vu{p}{1}$ we first generate the vector
\begin{equation}
\Vu{u}{1} = \V{x} \circ \Vu{p}{0},
\end{equation}
which we project onto the orthogonal complement of $\Vu{p}{0}$ to meet the orthogonality condition, yielding 
\begin{equation}
\Vu{\hat{p}}{1} = \left(\M{I} - \Vu{p}{0} \Vu{p}{0}^\transpose \M{W}_{\V{y}}\right) \Vu{u}{1}.
\end{equation}
This is a scaled version of $\Vu{p}{1}$.
To get the slope for the derivatives of the basis functions right, the vector $\V{x}'$ is calculated as
\begin{equation}
 \V{x}' = \sqrt{\V{e}^\transpose\M{W}_{\V{y}}\V{e} } \left(\frac{\Su{\hat{p}}{1,n} - \Su{\hat{p}}{1,1}}{\Su{x}{n} - \Su{x}{1}}\right)\V{e}.
\end{equation}
Using this, we calculate a scaled version of $\Vu{p}{1}'$ based on \Eq\eqref{eq:SubsVecDeriv} and \eqref{eq:SubsFullReccurenceVecDeriv}, yielding
\begin{equation}
\Vu{\hat{p}}{1}' =  \V{x}' \circ \Vu{p}{0}.
\end{equation}
To generate the basis functions fulfilling the norm conditions, we calculate the scaling factor
\begin{equation}
\Su{\alpha}{1} = \sqrt{\frac{1}{\Vu{\hat{p}}{1}^\transpose  \M{W}_{\V{y}}\Vu{\hat{p}}{1} + \Vu{\hat{p}}{1}'^\transpose\M{W}_{\md \V{y}} \Vu{\hat{p}}{1}'}}.
\end{equation}
From this we calculate the second pair of basis functions
\begin{equation}
\Vu{p}{1} =  \Su{\alpha}{1} \Vu{\hat{p}}{1}
\end{equation}
and
\begin{equation}
\Vu{p}{1}' = \Su{\alpha}{1} \Vu{\hat{p}}{1}'.
\end{equation}
This is the prerequisite to start the synthesis of higher order basis functions. The final set of basis functions  $\M{P} \defas\Mu{P}{d}$ and $\M{P}' \defas \Mu{P}{d}'$ of a certain degree $d$ can now be used in \Eq\eqref{eq:alpha}, \eqref{eq:valueApprox1} and \eqref{eq:derivApprox1} to calculate the coefficients and to finally approximate perturbed values and its derivatives with a polynomial of degree $d$.
A set of basis functions is shown in \Fig\ref{fig:BasisFns}.

The generation of a basis function set using the presented method is a generalized Gram-Schmidt process. \todo{cite GramSchmidt?}
As it can be seen, the generation of the basis depends only on the relative locations of the $x$~values\footnote{To improve numerical stability, $\V{x}$ is transformed to be centered at the origin and scaled to unit norm.} and not on the values itself.
If the abscissa values $\V{x}$ do not change, the basis function set can be calculated a priori and the solution to the weighted fitting problem reduces to a simple matrix-vector multiplication (see \Eq\eqref{eq:alpha}).
This is advantageous when implemented in smart sensors or low power controllers.
A further advantage is, that the weighting matrices $\M{W}_{\V{y}}$ and $\M{W}_{\md\V{y}}$ can be rank-deficient, e.g., points can be weighted with $0$ if they should not be considered. This can be helpful to suppress outliers.
\subsection{Covariance Propagation}
\label{ssec:CovProp}
Since we are dealing with noisy data, covariance propagation is a prerequisite for making assumptions about the quality of the approximated values.
Based on \Eq\eqref{eq:alpha} the covariance $\Mu{\Lambda}{\V{\gamma}}$ for the coefficients $\V{\gamma}$ is calculated as
\begin{equation}
	\Mu{\Lambda}{\V{\gamma}} =\M{P}_c^\transpose \M{W}_{c}  \Mu{\Lambda}{c} \M{W}_{c}^\transpose \M{P}_c
\end{equation}
with the block matrices
\begin{IEEEeqnarray}{rClrClrCl}
	\Mu{\Lambda}{c} &=&
	\begin{bmatrix}
		\M{\Lambda}_{\V{y}} & \M{0} \\ \M{0} & \M{\Lambda}_{\md \V{y}}
	\end{bmatrix}
	,& \ 
	\Mu{W}{c} &=&
	\begin{bmatrix}
		\M{W}_{\V{y}} & \M{0} \\ \M{0} & \M{W}_{\md \V{y}}
	\end{bmatrix},
	& \ 
	\Mu{P}{c} &=&
	\begin{bmatrix}
		\M{P} \\ \M{P}'
	\end{bmatrix}.
\end{IEEEeqnarray}

Similarly, the covariance matrices for the approximated values and derivatives can be propagated as \todo{difference beween lamda y and lamda y}
\begin{equation}
\Mu{\Lambda}{\Vt{y}} = \M{P} \M{P}_c^\transpose \M{W}_{c}  \Mu{\Lambda}{c} \M{W}_{c}^\transpose \M{P}_c \M{P}^\transpose
\end{equation}
and
\begin{equation}
\Mu{\Lambda}{\md \Vt{y}} = \M{P}' \M{P}_c^\transpose \M{W}_{c}  \Mu{\Lambda}{c} \M{W}_{c}^\transpose \M{P}_c \M{P}'^\transpose.
\end{equation}
Since $\Mu{W}{c}\defas\Mu{\Lambda}{c}^{-1}$  the middle-term results in $\M{W}_{c}  \Mu{\Lambda}{c} \M{W}_{c}^\transpose =   \M{W}_{c}$. Together with   the identity from \Eq\eqref{eq:BIdentity}, the above equations simplify to
\begin{equation}
\label{equ:covGamma}
\Mu{\Lambda}{\V{\gamma}} =\M{I},
\end{equation}
\begin{equation}
\Mu{\Lambda}{\Vt{y}} = \M{P} \M{P}^\transpose
\end{equation}
and
\begin{equation}
\Mu{\Lambda}{\md \Vt{y}} = \M{P}' \M{P}'^\transpose.
\end{equation}
These covariance matrices for the approximated coefficients and values, can be used to calculate confidence intervals or predictions intervals.
Note: As it can be seen in \Eq\eqref{equ:covGamma}, the presented method decorrelates the noise to i.i.d.~noise, to accord with Gauss's theorem.
%

%
\section{Numerical Example}
\label{sec:NumExample}
\begin{figure*}[!ht]
	\centering
	\includegraphics[width=\textwidth]{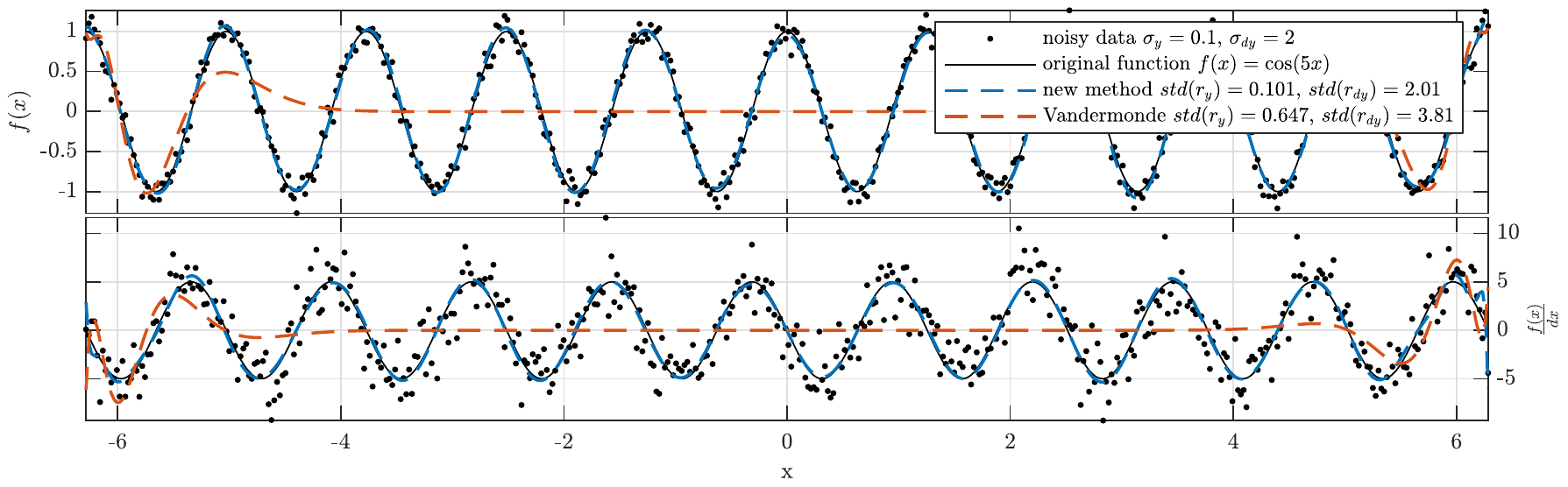}
	\caption{Approximation of a polynomial of degree $d=35$ to synthetic data generated from $f\left(x\right) = \cos\left(5x\right)$ with $\sigma_{y_i} = \sigma_{y} =  0.1$ and $\sigma_{\md y_i} = \sigma_{\md y} =  2$}
	\label{fig:CosFitting}
\end{figure*}

To test the validity of the herein presented method to approximate a polynomial given perturbed values and derivatives, a synthetic dataset is generated. The underlying function is defined as 
\begin{equation}
	f\left(x\right) = \cos\left(5x\right)
\end{equation}
with its analytical first derivative
\begin{equation}
\dx{f\left(x\right)}{x} = -5\sin\left(5x\right).
\end{equation}
The function and its derivative are evaluated in the range $\left[-2\pi, 2\pi\right]$ at $n=500$ equally spaced nodes, yielding the vectors of values and derivatives $\V{y}$ and $\V{y}'$.
Gaussian noise with different gains $\sigma_{y_i} = \sigma_{y} =  0.1$ and $\sigma_{\md y_i} = \sigma_{\md y} =  2$  is added to those vectors yielding the noisy measurement vectors $\Vh{y}$ and $\Vh{y}'$. \todo{should i use covariance matrix instead of $\sigma_{y}$ }
A polynomial of degree $d=35$ is used for approximating the noisy data set.
To test the developed method, a Monte~Carlo simulation is performed with $n_{\text{iter}} = 1000$ iterations.
As a measure, the standard deviation of the residuals $\std\{\Vu{r}{y}\}= \std\{\V{y} - \Vh{y}\}$ and  $\std\{\Vu{r}{\md y}\} = \std\{\V{y}' - \Vh{y}'\}$ are calculated in each run. Since the presented method uses covariance weighting, the standard deviation of the result should be the same as $\sigma_{y} \approx \std\{\V{y} - \Vh{y}\}$ and $\sigma_{\md y}\approx \std\{\V{y}' - \Vh{y}'\}$.
The mean value of the standard deviations over all runs is shown in \Fig\ref{fig:CosFitting}.
As it can be seen, although the noise gains are very different, the  presented method which uses covariance weighted approximation delivers the correct results for both, values and derivatives demonstrating the method to be valid.

As expected, the method using the Vandermonde basis as presented in \cite{Ritt2019}, is not stable for such a high degree.
As one can inspect, the approximation does not follow the signal. 
\section{Numerical Quality of Basis}
\label{ssec:NumQuality}
\begin{figure}[!ht]
	\centering
	\includegraphics{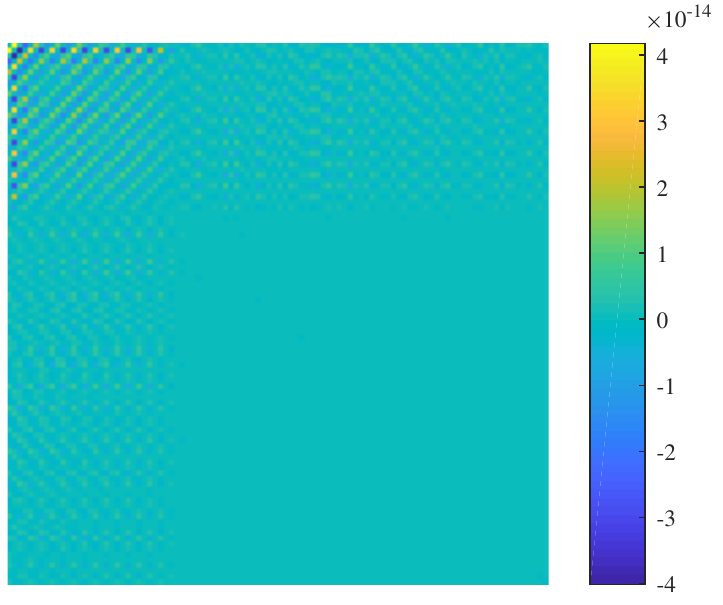}
	\caption{The structure of a residual matrix $\M{R} = \M{I} - \M{P}_c^\transpose \M{W}_{c} \M{P}_c$ for $d=100$ and $n=50$ with $\sigma_{y_i} = \sigma_{y} =  0.2$ and $\sigma_{\md y_i} = \sigma_{\md y} =  0.8$ }
\label{fig:numericalErrorPic}
\end{figure}

To verify the numerical quality of herein presented method, a meaningful measure has to be found.
As Wilkinson~\cite{Wilkinson1971} pointed out, a posteriori estimation of error bounds is preferred to a~priori error predictions in such cases.
The identity in \Eq\eqref{eq:BIdentity} can be written in terms of the block matrices as
\begin{equation}
	\M{P}_c^\transpose \M{W}_{c} \M{P}_c = \M{I}.
\end{equation}
Rewriting in terms of a unitary matrix $\M{U}$ yields
\begin{equation}
\M{U}^\transpose \M{U} = \M{I},
\end{equation}
with
\begin{equation}
\M{U} = \M{W}_{c}^{\frac{1}{2}} \M{P}_c.
\end{equation}
The residual matrix
\begin{equation}
	\M{R} = \M{I} - \M{P}_c^\transpose \M{W}_{c} \M{P}_c \approx \M{0}
\end{equation}
should contain zeros in a perfect reconstruction.
Due to numerical errors, this does not hold (see \Fig\ref{fig:numericalErrorPic}).
To summarize the numerical quality of the generated basis, the following error measures are tested in order to find the appropriate measure:
\begin{enumerate}
	\item \textbf{Maximum norm}. The maximum norm is the largest single element within the residual matrix, i.e.,  $\epsilon_{\text{max}} = \norm{\M{R}}_\text{max} = \max\{ |r_{ij}|\}$. Since 	this norm depends only on one specific entry, this may lead to wrong conclusions.
	\item \textbf{Frobenius norm}. The Frobenius norm is the square root of the sum of the squares of all entries in the residual matrix, i.e.,  $\epsilon_{\text{F}} = \norm{\M{R}}_\text{F} = \sqrt{\sum_i \sum_j r_{ij}^2}$. This norm is a measure for the total error.
	\item \textbf{Determinant}. The determinant of a matrix is a theoretical quality measure and equals $1$ for an ideal unitary matrix. The variation from this, i.e., $\epsilon_{\text{det}} = 1 - \det \{U\}$ is a measure for the quality of the tested basis function set.
	\item \textbf{Condition number}. The condition number of a matrix is connected to the error propagation. For a unitary matrix the condition number is $1$. Thus, the measure tested is $\epsilon_{\text{cond}} = 1 - \cond \{U\}$.
	\item \textbf{Rank}. Since the presented method generates a weighted orthogonal basis function set, the rank of $\M{U}$ should be full rank, i.e., no linear dependencies. Due to round off errors this can be used to find linear dependencies, i.e., $\epsilon_{\text{rank}} = n - \rank{\M{U}}$ \todo{maybe remove}.
\end{enumerate}

\begin{figure}[!ht]
	\centering
	\includegraphics{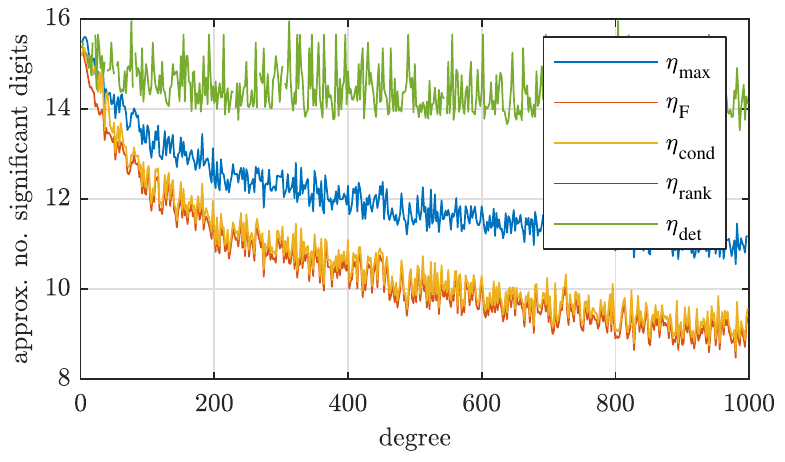}
	\caption{Comparison of different error measures for a complete basis function set with  $d = 2n - 1$ and $\sigma_{y_i} = \sigma_{y} =  0.2$ and $\sigma_{\md y_i} = \sigma_{\md y} =  0.8$}
	\label{fig:DifferentQualityMeasures}
\end{figure}

Based on those measures the approximate number of significant digits $\eta$ is calculated as
\begin{equation}
	\eta_m = -\log_{10} \left(\epsilon_m\right).
\end{equation}
These measures are calculated for a complete basis function set where  the number of basis functions equals the number of data points (values and derivatives). Therefore, the degree of the resulting polynomial is $d = 2n - 1$. In \Fig\ref{fig:DifferentQualityMeasures} the different measures are presented for varying degrees.
 As it can be seen, the Frobenius~norm and the condition number are the most meaningful measures, since they show the highest dependency on the degree of the resulting polynomial. The rank measure is not visible, since the proposed algorithm generates full-rank basis function sets, so there is no error visible. \todo{maybe remove} 
 Since the Frobenius~norm is a measure for the total error, this measure is chosen to compare to other algorithms in the following.
 
\begin{figure}[!ht]
 	\centering
 	\includegraphics{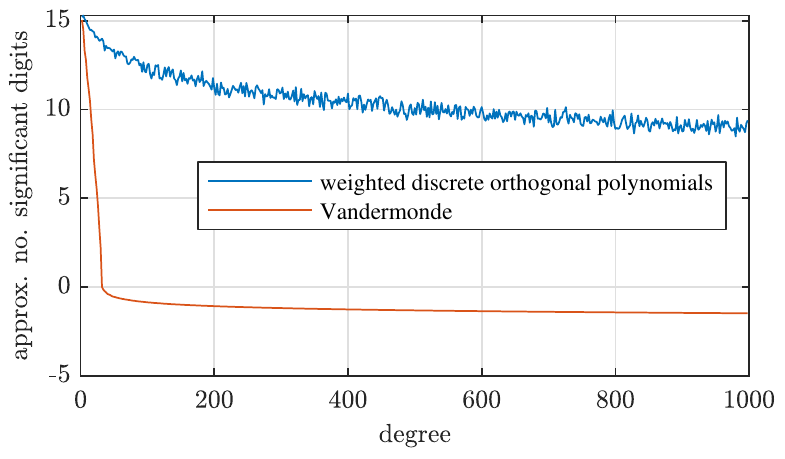}
 	\caption{Numerical quality of Vandermonde basis compared to the new method using a complete basis with $\sigma_{y_i} = \sigma_{y} =  0.2$ and $\sigma_{\md y_i} = \sigma_{\md y} =  0.8$}
 	\label{fig:FullBasisQuality}
 \end{figure}

 The new method is compared to the one presented in \cite{Ritt2019}. It uses a Vandermonde basis function set to solve the same problem. As \cite{Harker2013a} pointed out, the Vandermonde basis for normal polynomial regression gets degenerate at high degrees. This behaviour can also be inspected within this paper. As it can be seen in \Fig\ref{fig:FullBasisQuality}, the new method generates a more stable result also for high degrees.

\begin{figure}[!ht]
	\centering
	\includegraphics{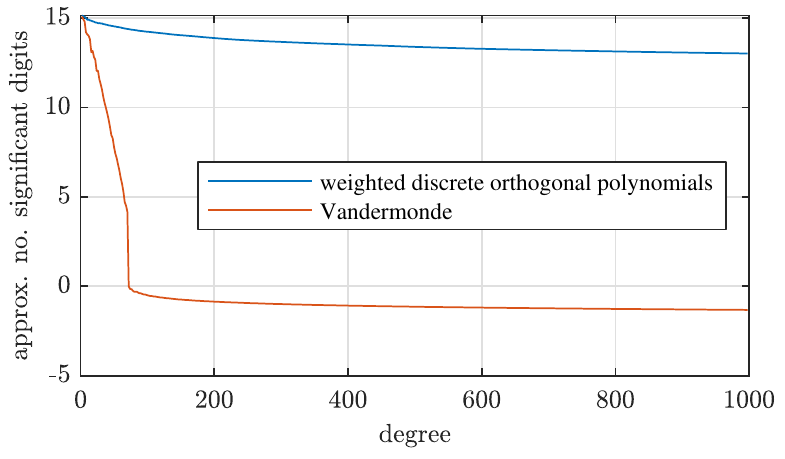}
	\caption{Numerical quality of Vandermonde basis compared to the new method using an incomplete  basis with $\sigma_{y_i} = \sigma_{y} =  0.2$, $\sigma_{\md y_i} = \sigma_{\md y} =  0.8$ and $n=1000$}
	\label{fig:IncompleteBasisQuality}
\end{figure}

 Since this new method can also be used for approximation (overdetermined system of equations), the quality of the basis is determined for an incomplete basis function set for $n=1000$.
 The result for various degrees of polynomial is visualized in \Fig\ref{fig:IncompleteBasisQuality}, showing that the new method performs also better for an incomplete basis.
\section{Conclusion}
\label{sec:Conclusion}
The herein presented method introduces a novel polynomial fitting framework for the approximation of value and derivative data.
Including both sources of information within the fitting procedure improves the quality of the fitted polynomial improving both, reconstruction of values and derivatives. 
The method uses a recurrence relation with full re-orthogonalization together with covariance weighting for introducing a metric between value and derivative domain, yielding a set of discrete orthogonal polynomials.
As it is shown, the generated basis function set is numerically more stable compared to other methods, especially for high degree polynomials.
Using this basis, the fitting problem is reduced to inner products, which is beneficial in terms of computational efficiency.
Due to the covariance weighting, the noise associated with the channels is decorrelated to i.i.d.~noise, to accord with Gauss's theorem.
Thus, there is no bias based on different noise parameters.
The validity of the method is tested and presented on a numerical example, where a polynomial of degree $d=35$ is fitted to a periodic function.
In future research, this method will be adapted for fitting data given noisy constraints, based on discrete orthogonal polynomials.
\section*{Acknowledgement}
This work was partially funded under the auspices of the EIT~-~KIC Raw materials program within the project “Maintained Mine and Machine”
(MaMMa) with the grant agreement number: [EIT/RAW MATERIALS/SGA2018/1].
\bibliographystyle{IEEEtran}
\bibliography{LiteraturPublicatons-2018-PaperHermiteApproxDOP}
\end{document}